\documentclass[11pt,oneside,leqno]{article}
\usepackage[utf8]{inputenc}
\usepackage[russian]{babel}
\usepackage{amssymb,amsmath,amsthm}
\renewcommand{\subsection}{\refstepcounter{subsection}%
\par\bigskip\noindent\textbf{\upshape\thesubsection. }}
\renewcommand{\subsubsection}{\refstepcounter{subsubsection}%
\par\medskip\noindent\textbf{\upshape\thesubsubsection.  }}
\renewcommand{\paragraph}{\refstepcounter{paragraph}%
\par\smallskip\noindent\textbf{\upshape\theparagraph. }}

\numberwithin{equation}{subsection}

\renewcommand{\thesubsection}{\arabic{subsection}}

\newcommand{\Wo}{{\raisebox{0.2ex}{\(\stackrel{\circ}{W}\)}}{}}
\newcommand{\ind}{\operatorname{ind}}

\title{\begin{flushleft}\normalsize УДК~517.984\end{flushleft}
Асимптотика собственных значений задачи Штурма-Лиувилля с дискретным
самоподобным весом}
\author{А.~А.~Владимиров, И.~А.~Шейпак\footnote{%
Работа поддержана РФФИ, грант \No~07-01-00283, фондом поддержки
ведущих научных школ, грант~НШ-5247.2006.1, и фондом INTAS,
грант~05-1000008-7883.}}
\date{}
\begin{document}
\renewcommand{\proofname}{{\upshape Д\,о\,к\,а\,з\,а\,т\,е\,л\,ь\,с\,т\,в\,о.}}
\maketitle
\begin{abstract}
В статье изучается вопрос об асимптотике спектра граничной задачи
\begin{gather*}
	-y''-\lambda\rho y=0,\\
	y(0)=y(1)=0
\end{gather*}
в случае, когда вес \(\rho\in\Wo_2^{-1}[0,1]\) представляет собой обобщённую
производную самоподобной функции \(P\in L_2[0,1]\) нулевого спектрального порядка.
\end{abstract}

\section{Введение}\label{par:1}
\subsection\label{pt:1}
В настоящей статье нами будет продолжено начатое в работах \cite{VSh1} и \cite{VSh2}
изучение спектральных асимптотик граничной задачи
\begin{gather}\label{eq:1}
	-y''-\lambda\rho y=0,\\ \label{eq:2}
	y(0)=y(1)=0,
\end{gather}
где вес \(\rho\in\Wo_2^{-1}[0,1]\) представляет собой обобщённую производную
самоподобной функции \(P\in L_2[0,1]\). Говоря более точно, мы намерены теперь
подвергнуть исследованию ранее не рассмотренный случай, когда спектральный порядок
\(D\) функции \(P\) равен нулю. Определение понятия спектрального порядка
квадратично суммируемой самоподобной функции было дано в \cite[\S~3.3]{VSh1}.

\subsection
Особенностью случая \(D=0\) является неприменимость использованной в работах
\cite{SV}, \cite{VSh1} и \cite{VSh2} при рассмотрении случая \(D>0\) техники,
основанной на теории восстановления. Действительно, указанная техника существенно
опирается на факт экспоненциального убывания некоторых вспомогательных функций
(см., например, оценку \cite[(4.11)]{VSh1}), между тем как в случае \(D=0\)
для таких функций имеется лишь неулучшаемая, вообще говоря, оценка порядка
\(O(1)\). В соответствии со сказанным, далее нами будет развита другая методика
исследования.

\subsection
В дальнейшем через \(\mathfrak H\) мы будем обозначать пространство Соболева
\(\Wo_2^1[0,1]\), снабжённое скалярным произведением
\[
	\langle y,z\rangle\rightleftharpoons\int\limits_0^1 y'\overline{z'}\,dx.
\]
Через \(\mathfrak H'\) мы при этом будем обозначать пространство, двойственное
к \(\mathfrak H\) относительно \(L_2[0,1]\), т.~е. получаемое пополнением
пространства \(L_2[0,1]\) по норме
\[
	\|y\|_{\mathfrak H'}\rightleftharpoons\sup\limits_{\|z\|_{\mathfrak H}=1}
	\left|\int\limits_0^1 y\overline{z}\,dx\right|.
\]
Если рассмотреть оператор вложения \(J:\mathfrak H\to L_2[0,1]\), то непосредственно
из определения пространства \(\mathfrak H'\) вытекает возможность непрерывного
продолжения сопряжённого оператора \(J^*:L_2[0,1]\to\mathfrak H\) до изометрии
\(J^+:\mathfrak H'\to\mathfrak H\).

Как и в предшествующих работах \cite{VSh1} и \cite{VSh2}, в качестве операторной
модели задачи \ref{pt:1}\,\eqref{eq:1},~\ref{pt:1}\,\eqref{eq:2} мы будем рассматривать
линейный пучок \(T_{\rho}:\mathfrak H\to\mathfrak H'\) ограниченных операторов,
удовлетворяющий тождеству
\[
	(\forall\lambda\in\mathbb R)\,(\forall y\in\mathfrak H)\qquad
	\langle J^+T_{\rho}(\lambda)y,y\rangle=\int\limits_0^1\left(|y'|^2+
	\lambda P\cdot(|y|^2)'\right)\,dx.
\]

\subsection
Через \(\ind E\) мы далее будем обозначать отрицательный индекс инерции действующего
в некотором гильбертовом пространстве \(\mathfrak E\) ограниченного эрмитова
оператора \(E\), т.~е. точную верхнюю грань размерностей подпространств
\(\mathfrak M\subseteq\mathfrak E\), удовлетворяющих условию
\[
	(\exists\varepsilon>0)\,(\forall y\in\mathfrak M)\qquad
	\langle Ey,y\rangle\leqslant-\varepsilon\,\|y\|^2_{\mathfrak E}.
\]

Обозначения \(n\), \(a_k\), \(d_k\) и \(\beta_k\), где \(k=1,\ldots, n\), мы будем
использовать для записи параметров самоподобия функции \(P\). Определение этих
параметров и описание накладываемых на них условий дано в \cite[\S~3.2]{VSh1}.
Особо отметим, что указанные условия гарантируют выполнение неравенств
\(a_k|d_k|<1\) (ср.~далее формулировки утверждений из параграфа~\ref{par:3}).

Через \(\alpha_k\), где \(k=0,1,\ldots, n\), мы будем обозначать величины,
определяемые рекуррентными соотношениями \(\alpha_0\rightleftharpoons 0\)
и \(\alpha_k\rightleftharpoons\alpha_{k-1}+a_k\) при \(k\neq 0\).

Наконец, через \(m\) мы будем обозначать натуральное число \(m\in [1,n]\),
удовлетворяющее условию
\[
	(\forall k\in [1,n]\setminus\{m\})\qquad d_k=0
\]
(ср. \cite[\S~3.3]{VSh1}). При этом, ввиду тривиальности случая \(d_m=0\),
будет обычно предполагаться выполнение неравенства \(d_m\neq 0\).

\subsection
В дальнейшем при ссылках на разделы статьи, не принадлежащие параграфу,
внутри которого даётся ссылка, будет дополнительно указываться номер
параграфа. При ссылках на формулы, не принадлежащие пункту, внутри
которого даётся ссылка, будет дополнительно указываться номер пункта.

\section{Вспомогательные понятия и утверждения}\label{par:2}
\subsection
Рассмотрим два подпространства \(\mathfrak H_1\subseteq\mathfrak H\)
и \(\mathfrak H_2\subseteq\mathfrak H\), определяемые следующим образом.
Подпространство \(\mathfrak H_1\) имеет вид
\[
	\mathfrak H_1\rightleftharpoons\{y\in\mathfrak H\mid
	(\forall x\not\in [\alpha_{m-1},\alpha_m])\quad y(x)=0\}.
\]
Подпространство \(\mathfrak H_2\) представляет собой \((n-1)\)-мерную линейную
оболочку функций \(e_k\in\mathfrak H\), где \(k=1,\ldots, n-1\), имеющих вид
\[
	e_k(x)=\left\{
		\begin{aligned}
			&\dfrac{x-\gamma_k}{\alpha_k-\gamma_k}&
			&\text{при }x\in [\gamma_k,\alpha_k],\\
			&\dfrac{\delta_k-x}{\delta_k-\alpha_k}&
			&\text{при }x\in [\alpha_k,\delta_k],\\
			&0 &&\text{иначе,}
		\end{aligned}
	\right.
\]
где положено
\begin{align*}
	\gamma_k&\rightleftharpoons\left\{
		\begin{aligned}
			&\alpha_{k-1}&&\text{при }k\neq m,\\
			&\alpha_m-a_ma_n&&\text{при }k=m,
		\end{aligned}
	\right.\\
	\delta_k&\rightleftharpoons\left\{
		\begin{aligned}
			&\alpha_{k+1}&&\text{при }k\neq m-1,\\
			&\alpha_{m-1}+a_ma_1&&\text{при }k=m-1.
		\end{aligned}
	\right.
\end{align*}
Имеют место следующие два факта:

\subsubsection\label{2:1}
{\itshape Ортогональное дополнение прямой суммы \(\mathfrak H_1\dotplus
\mathfrak H_2\) допускает представление в виде
\[
	\mathfrak H\ominus(\mathfrak H_1\dotplus\mathfrak H_2)=
	\left\{y\in\mathfrak H\;\vrule\;
	\left(\forall x\in [\alpha_{m-1},\alpha_m]\cup
	\bigcup\limits_{k=1}^{n-1} \{\alpha_k\}\right)\quad y(x)=0\right\}.
\]
}

\bigskip
Справедливость утверждения \ref{2:1} легко устанавливается прямым вычислением.

\subsubsection\label{2:2}
{\itshape Пусть \(\lambda\) "--- вещественное число, а \(\mathfrak M\subseteq
\mathfrak H\) "--- конечномерное подпространство, на котором квадратичная форма
оператора \(J^+T_{\rho}(\lambda)\) отрицательна. Тогда существует подпространство
\(\mathfrak N\subseteq\mathfrak H_1\dotplus\mathfrak H_2\) размерности
\(\dim\mathfrak M\), на котором квадратичная форма оператора
\(J^+T_{\rho}(\lambda)\) также отрицательна.
}

\begin{proof}
Доказываемое утверждение с очевидностью вытекает из тождества
\[
	(\forall y\in\mathfrak H_1\dotplus\mathfrak H_2)\,
	(\forall z\perp\mathfrak H_1\dotplus\mathfrak H_2)\qquad
	\langle J^+T_{\rho}(\lambda)\,(y+z),(y+z)\rangle=
	\langle J^+T_{\rho}(\lambda)y,y\rangle+\|z\|^2_{\mathfrak H},
\]
легко получаемого на основе утверждения \ref{2:1} и факта самоподобия функции \(P\).
\end{proof}

\subsection
Рассмотрим два линейных пучка \(A:\mathfrak H_1\to\mathfrak H_1\)
и \(C:\mathfrak H_2\to\mathfrak H_2\) ограниченных операторов, удовлетворяющие
тождествам
\begin{align*}
	(\forall\lambda\in\mathbb R)\,(\forall y\in\mathfrak H_1)\qquad
	\langle A(\lambda)y,y\rangle&=\int\limits_0^1\left(|y'|^2+\lambda P\cdot
	(|y|^2)'\right)\,dx,\\
	(\forall\lambda\in\mathbb R)\,(\forall y\in\mathfrak H_2)\qquad
	\langle C(\lambda)y,y\rangle&=\int\limits_0^1\left(|y'|^2+\lambda P\cdot
	(|y|^2)'\right)\,dx,
\end{align*}
а также оператор \(B:\mathfrak H_1\to\mathfrak H_2\), удовлетворяющий тождеству
\[
	(\forall y\in\mathfrak H_1)\,(\forall z\in\mathfrak H_2)\qquad
	\langle By,z\rangle=\int\limits_0^1 y'\overline{z'}\,dx.
\]
Имеют место следующие два факта:

\subsubsection\label{3:1}
{\itshape Пусть \(\lambda\) "--- вещественное число. Тогда выполняется равенство
\[
	\ind A(\lambda)=\ind J^+T_{\rho}(a_md_m\,\lambda).
\]
}

\begin{proof}
Рассмотрим изометрический оператор \(U:\mathfrak H_1\to\mathfrak H\), определяемый
тождеством
\[
	(\forall y\in\mathfrak H_1)\,(\forall x\in [0,1])\qquad
	[U y](x)=\dfrac{y(\alpha_{m-1}+a_mx)}{\sqrt{a_m}}.
\]
Доказываемое утверждение с очевидностью вытекает из тождества
\[
	(\forall y\in\mathfrak H_1)\qquad \langle A(\lambda)y,y\rangle=
	\langle J^+T_{\rho}(a_md_m\,\lambda)U y,U y\rangle,
\]
легко устанавливаемого непосредственными выкладками с учётом самоподобия
функции \(P\).
\end{proof}

\subsubsection\label{3:2}
{\itshape Пусть \(\lambda\) "--- вещественное число, не принадлежащее спектру пучка
\(C\). Тогда выполняется равенство
\[
	\ind J^+T_{\rho}(\lambda)=\ind [A(\lambda)-B^*C^{-1}(\lambda)B]+
	\ind C(\lambda).
\]
}

\begin{proof}
Прямым вычислением легко устанавливается, что для любых функций
\(y\in\mathfrak H_1\) и \(z\in\mathfrak H_2\) выполняется равенство
\[
	\langle J^+T_{\rho}(\lambda)\,(y+z),(y+z)\rangle=
	\langle [A(\lambda)-B^*C^{-1}(\lambda)B]y,y\rangle+
	\langle C(\lambda)u,u\rangle,
\]
где положено \(u\rightleftharpoons z+C^{-1}(\lambda)By\). Доказываемое утверждение
вытекает из этого равенства и утверждения \ref{2:2}.
\end{proof}

\subsection
Рассмотрим величины \(\zeta_k\), где \(k=1,\ldots, n-1\), имеющие вид
\[
	\zeta_k\rightleftharpoons\left\{
		\begin{aligned}
			&\beta_m-\beta_{m-1}+d_m\beta_1&&\text{при }k=m-1,\\
			&\beta_{m+1}-\beta_m-d_m\beta_n&&\text{при }k=m,\\
			&\beta_{k+1}-\beta_k&&\text{иначе.}
		\end{aligned}
	\right.
\]
Обозначим также через \(\mathrm Z_{\pm}\) две величины
\[
	\mathrm Z_{\pm}\rightleftharpoons\#\{k\in [1,n-1]\mid\pm\zeta_k>0\}.
\]
Имеют место следующие два факта:

\subsubsection\label{4:1}
{\itshape Для любого достаточно большого вещественного числа \(\lambda>0\)
выполняется равенство
\[
	\ind C(\lambda)=\mathrm Z_+.
\]
}

Утверждение \ref{4:1} немедленно вытекает из легко проверяемого тождества
\begin{equation}\label{eq:3}
	(\forall\lambda\in\mathbb R)\,(\forall y\in\mathfrak H_2)\qquad
	\langle C(\lambda)y,y\rangle=\|y\|^2_{\mathfrak H}-
	\lambda\sum\limits_{k=1}^{n-1}\zeta_k\cdot|y(\alpha_k)|^2.
\end{equation}

\subsubsection\label{4:2}
{\itshape Пусть выполнено равенство \(\mathrm Z_++\mathrm Z_-=n-1\). Тогда для любого достаточно
большого вещественного числа \(\lambda>0\) оператор \(C(\lambda)\) является
ограниченно обратимым, причём при \(\lambda\to+\infty\) справедлива асимптотика
\[
	\|C^{-1}(\lambda)\|=O(\lambda^{-1}).
\]
}

Утверждение \ref{4:2} также представляет собой несложное следствие тождества
\eqref{eq:3}.

\section{Основные результаты}\label{par:3}
\subsection
Имеют место следующие три факта:

\subsubsection\label{tm:1}
{\itshape Пусть выполняются соотношения \(d_m>0\), \(\mathrm Z_+>0\) и \(\mathrm Z_++\mathrm Z_-=n-1\).
Тогда существуют вещественные числа \(\mu_l>0\), где \(l=1,\ldots, \mathrm Z_{+}\),
для которых последовательность \(\{\lambda_k\}_{k=1}^{\infty}\) занумерованных
в порядке возрастания положительных собственных значений задачи
\ref{par:1}.\ref{pt:1}\,\eqref{eq:1}, \ref{par:1}.\ref{pt:1}\,\eqref{eq:2}
удовлетворяет при \(k\to\infty\) асимптотикам
\[
	\lambda_{l+k\mathrm Z_{+}}=\mu_l\cdot (a_md_m)^{-k}\cdot (1+o(1)).
\]
}

\begin{proof}
Согласно утверждению \ref{par:2}.\ref{4:2}, найдётся такое вещественное число
\(\lambda_0>0\), что при любом \(\lambda>\lambda_0\) будет выполняться неравенство
\(\|C^{-1}(\lambda)\|<\lambda_0/(3\lambda)\). Отсюда и из очевидной оценки
\(\|B\|\leqslant 1\) следует, что при любом \(\lambda>\lambda_0\) будут также
выполняться неравенства
\begin{equation}\label{eq:4}
	\ind [A(\lambda)+\lambda_0/(3\lambda)]\leqslant
	\ind [A(\lambda)-B^*C^{-1}(\lambda)B]\leqslant
	\ind [A(\lambda)-\lambda_0/(3\lambda)].
\end{equation}
При этом с использованием очевидных равенств
\[
	\ind [A(\lambda)\pm\lambda_0/(3\lambda)]=
	\ind \bigl[(1\pm\lambda_0/(3\lambda))^{-1}\cdot
	(A(\lambda)\pm\lambda_0/(3\lambda))\bigr]
\]
из оценок \eqref{eq:4} легко выводятся оценки
\begin{equation}\label{eq:5}
	\ind A(\lambda-\lambda_0/2)\leqslant
	\ind [A(\lambda)-B^*C^{-1}(\lambda)B]\leqslant
	\ind A(\lambda+\lambda_0/2).
\end{equation}

Рассмотрим теперь неубывающую функцию \(s:\mathbb R\to\mathbb R\) вида
\[
	(\forall t\in\mathbb R)\qquad s(t)=\#\{k\geqslant 1\mid\lambda_k<e^t\}.
\]
Согласно теореме \cite[Теорема~4.1]{VSh1}, эта функция удовлетворяет тождеству
\[
	(\forall t\in\mathbb R)\qquad s(t)=\ind J^+T_{\rho}(e^t).
\]
Объединяя последнее тождество с оценками \eqref{eq:5} и утверждениями
\ref{par:2}.\ref{3:1}, \ref{par:2}.\ref{3:2} и \ref{par:2}.\ref{4:1},
устанавливаем, что при любом \(t>\ln\lambda_0\) будут справедливы оценки
\begin{equation}\label{eq:6}
	s(t+\ln(a_md_m)-\lambda_0 e^{-t})+\mathrm Z_+\leqslant s(t)\leqslant
	s(t+\ln(a_md_m)+\lambda_0 e^{-t})+\mathrm Z_+.
\end{equation}

Из оценок \eqref{eq:6} следует, что найдётся точка \(t_0>\ln\lambda_0\), для которой
окрестность радиуса \(2\lambda_0 e^{-t_0}/(1-a_md_m)\) не содержит точек разрыва
функции \(s\). Объединяя те же оценки с фактом простоты всех собственных значений
задачи \ref{par:1}.\ref{pt:1}\,\eqref{eq:1}, \ref{par:1}.\ref{pt:1}\,\eqref{eq:2}
(см., например, теорему \cite[Теорема~4.1]{VSh1}), убеждаемся, что при любом
\(k\in\mathbb N\) отрезок
\[
	\Delta_k\rightleftharpoons [t_0-k\cdot\ln(a_md_m),
	t_0-(k+1)\cdot\ln(a_md_m)]
\]
содержит ровно \(\mathrm Z_+\) точек разрыва функции \(s\), причём каждая из них
является внутренней точкой отрезка \(\Delta_k\). Кроме того, при любых \(k\gg 0\)
и \(l=1,\ldots,\mathrm Z_+\) расстояние между \(l\)-тыми слева точками разрыва функции
\(s\) на отрезках \(\Delta_k\) и \(\Delta_{k+1}\) отличается от величины
\(-\ln(a_md_m)\) не более, чем на \(\lambda_0 e^{-t_0}\cdot (a_md_m)^{k+1}\).
Тем самым, доказываемое утверждение справедливо.
\end{proof}

\subsubsection\label{tm:2}
{\itshape Пусть выполняются соотношения \(d_m>0\), \(\mathrm Z_->0\) и \(\mathrm Z_++\mathrm Z_-=n-1\).
Тогда существуют вещественные числа \(\mu_l>0\), где \(l=1,\ldots, \mathrm Z_{-}\),
для которых последовательность \(\{\lambda_{-k}\}_{k=1}^{\infty}\) занумерованных
в порядке убывания отрицательных собственных значений задачи
\ref{par:1}.\ref{pt:1}\,\eqref{eq:1}, \ref{par:1}.\ref{pt:1}\,\eqref{eq:2}
удовлетворяет при \(k\to\infty\) асимптотикам
\[
	\lambda_{-(l+k\mathrm Z_{-})}=-\mu_l\cdot (a_md_m)^{-k}\cdot (1+o(1)).
\]
}

Доказательство утверждения \ref{tm:2} аналогично доказательству утверждения
\ref{tm:1}.

\subsubsection\label{tm:3}
{\itshape Пусть выполняются соотношения \(d_m<0\) и \(\mathrm Z_++\mathrm Z_-=n-1\). Тогда
существуют вещественные числа \(\mu_l>0\), где \(l=1,\ldots, n-1\),
для которых последовательность \(\{\lambda_k\}_{k=1}^{\infty}\) занумерованных
в порядке возрастания положительных собственных значений задачи
\ref{par:1}.\ref{pt:1}\,\eqref{eq:1}, \ref{par:1}.\ref{pt:1}\,\eqref{eq:2}
удовлетворяет при \(k\to\infty\) асимптотикам
\[
	\lambda_{l+k(n-1)}=\mu_l\cdot (a_m|d_m|)^{-2k}\cdot (1+o(1)),
\]
а последовательность \(\{\lambda_{-k}\}_{k=1}^{\infty}\) занумерованных в порядке
убывания отрицательных собственных значений задачи
\ref{par:1}.\ref{pt:1}\,\eqref{eq:1}, \ref{par:1}.\ref{pt:1}\,\eqref{eq:2}
удовлетворяет при \(k\to\infty\) асимптотикам
\[
	\lambda_{-(l+\mathrm Z_{-}+k(n-1))}=-\mu_l\cdot (a_m|d_m|)^{-2k-1}\cdot
	(1+o(1)).
\]
}

Доказательство утверждения \ref{tm:3} также аналогично доказательству утверждения
\ref{tm:1}.

\section{Примеры}
\subsection
\begin{table}[t]
\begin{center}
\begin{tabular}{|r|r|rrr|rrr|}
\hline
{\(l\)}&{\(k\)}&\multicolumn{3}{|c|}{\(\lambda_{l+2k}\)}&
\multicolumn{3}{|c|}{\(\lambda_{l+2k}/6^{k}\)}\\ \hline
1&0&\(4,93\cdot 10^0\)&\(\pm\)&\(1\%\)&4,9341&\(\pm\)&\(10^{-4}\)\\
2&0&\(1,36\cdot 10^1\)&\(\pm\)&\(1\%\)&13,6598&\(\pm\)&\(10^{-4}\)\\
1&1&\(4,94\cdot 10^1\)&\(\pm\)&\(1\%\)&8,2322&\(\pm\)&\(10^{-4}\)\\
2&1&\(8,85\cdot 10^1\)&\(\pm\)&\(1\%\)&14,7576&\(\pm\)&\(10^{-4}\)\\
1&2&\(2,96\cdot 10^2\)&\(\pm\)&\(1\%\)&8,2330&\(\pm\)&\(10^{-4}\)\\
2&2&\(5,31\cdot 10^2\)&\(\pm\)&\(1\%\)&14,7577&\(\pm\)&\(10^{-4}\)\\
1&3&\(1,78\cdot 10^3\)&\(\pm\)&\(1\%\)&8,2330&\(\pm\)&\(10^{-4}\)\\
2&3&\(3,19\cdot 10^3\)&\(\pm\)&\(1\%\)&14,7577&\(\pm\)&\(10^{-4}\)\\
\hline
\end{tabular}
\end{center}

\vspace{0.5cm}
\caption{Оценки первых собственных значений для случая \(n=3\), \(a_1=a_2=
a_3=1/3\), \(m=3\), \(d_3=1/2\), \(\beta_1=0\), \(\beta_2=2/3\), \(\beta_3=1\).}
\label{tab:1}
\end{table}
В таблице \ref{tab:1} представлены результаты численных расчётов для первых восьми
положительных собственных значений задачи Штурма--Лиувилля, весовой функцией в которой
выступает обобщённая производная квадратично суммируемой функции с параметрами
самоподобия \(n=3\), \(a_1=a_2=a_3=1/3\), \(m=3\), \(d_3=1/2\), \(\beta_1=0\),
\(\beta_2=2/3\), \(\beta_3=1\). В этом случае выполняются равенства \(\zeta_1=2/3\),
\(\zeta_2=1/3\), \(\mathrm Z_+=2\), \(\mathrm Z_-=0\). Данные таблицы иллюстрируют утверждение
\ref{par:3}.\ref{tm:1}.

\subsection
\begin{table}[t]
\begin{center}
\begin{tabular}{|r|r|rrr|rrr|}
\hline
{\(l\)}&{\(k\)}&\multicolumn{3}{|c|}{\(-\lambda_{-(l+k)}\)}&
\multicolumn{3}{|c|}{\(-\lambda_{-(l+k)}/6^{k}\)}\\ \hline
1&0&\(5,10\cdot 10^0\)&\(\pm\)&\(1\%\)&5,1005&\(\pm\)&\(10^{-4}\)\\
1&1&\(2,60\cdot 10^1\)&\(\pm\)&\(1\%\)&4,3459&\(\pm\)&\(10^{-4}\)\\
1&2&\(1,56\cdot 10^2\)&\(\pm\)&\(1\%\)&4,3458&\(\pm\)&\(10^{-4}\)\\
1&3&\(9,39\cdot 10^2\)&\(\pm\)&\(1\%\)&4,3458&\(\pm\)&\(10^{-4}\)\\
\hline
\end{tabular}
\end{center}

\vspace{0.5cm}
\caption{Оценки первых собственных значений для случая \(n=3\), \(a_1=a_2=
a_3=1/3\), \(m=3\), \(d_3=1/2\), \(\beta_1=\beta_3=0\), \(\beta_2=-1\).}
\label{tab:2}
\end{table}
В таблице \ref{tab:2} представлены данные численных расчётов первых четырёх
отрицательных собственных значений задачи Штурма--Лиувилля, весовой функцией
в которой выступает обобщённая производная квадратично суммируемой функции
с параметрами самоподобия \(n=3\), \(a_1=a_2=a_3=1/3\), \(m=3\), \(d_3=1/2\),
\(\beta_1=\beta_3=0\), \(\beta_2=-1\). В этом случае выполняются равенства
\(\zeta_1=-1\), \(\zeta_2=1\), \(\mathrm Z_+=\mathrm Z_-=1\). Данные таблицы иллюстрируют
утверждение \ref{par:3}.\ref{tm:2}.

\subsection
\begin{table}[t]
\begin{center}
\begin{tabular}{|r|r|rrr|rrr|rrr|rrr|}
\hline
{\(l\)}&{\(k\)}&\multicolumn{3}{|c|}{\(\lambda_{l+2k}\)}&
\multicolumn{3}{|c|}{\(\lambda_{l+2k}/6^{2k}\)}&
\multicolumn{3}{|c|}{\(-\lambda_{-(l+1+2k)}\)}&
\multicolumn{3}{|c|}{\(-\lambda_{-(l+1+2k)}/6^{2k+1}\)}\\ \hline
1&0&\(4,31\cdot 10^0\)&\(\pm\)&\(1\%\)&4,3146&\(\pm\)&\(10^{-4}\)
&\(2,55\cdot 10^1\)&\(\pm\)&\(1\%\)&4,2572&\(\pm\)&\(10^{-4}\)\\
2&0&\(3,81\cdot 10^1\)&\(\pm\)&\(1\%\)&38,0536&\(\pm\)&\(10^{-4}\)
&\(2,28\cdot 10^2\)&\(\pm\)&\(1\%\)&38,0535&\(\pm\)&\(10^{-4}\)\\
1&1&\(1,53\cdot 10^2\)&\(\pm\)&\(1\%\)&4,2572&\(\pm\)&\(10^{-4}\)
&\(9,19\cdot 10^2\)&\(\pm\)&\(1\%\)&4,2572&\(\pm\)&\(10^{-4}\)\\
2&1&\(1,37\cdot 10^3\)&\(\pm\)&\(1\%\)&38,0535&\(\pm\)&\(10^{-4}\)
&\(8,22\cdot 10^3\)&\(\pm\)&\(1\%\)&38,0535&\(\pm\)&\(10^{-4}\)\\
1&2&\(5,52\cdot 10^3\)&\(\pm\)&\(1\%\)&4,2572&\(\pm\)&\(10^{-4}\)
&\(3,31\cdot 10^4\)&\(\pm\)&\(1\%\)&4,2572&\(\pm\)&\(10^{-4}\)\\
2&2&\(4,93\cdot 10^4\)&\(\pm\)&\(1\%\)&38,0535&\(\pm\)&\(10^{-4}\)
&\(2,96\cdot 10^5\)&\(\pm\)&\(1\%\)&38,0535&\(\pm\)&\(10^{-4}\)\\
\hline
\end{tabular}
\end{center}

\vspace{0.5cm}
\caption{Оценки первых собственных значений для случая \(n=3\), \(a_1=a_2=
a_3=1/3\), \(m=3\), \(d_3=-1/2\), \(\beta_1=\beta_3=0\), \(\beta_2=-1\).}
\label{tab:3}
\end{table}
В таблице \ref{tab:3} представлены данные численных расчётов первых шести
положительных и семи отрицательных (исключая первое) собственных значений задачи
Штурма--Лиувилля, весовой функцией в которой выступает обобщённая производная
квадратично суммируемой функции с параметрами самоподобия \(n=3\),
\(a_1=a_2=a_3=1/3\), \(m=3\), \(d_3=-1/2\), \(\beta_1=\beta_3=0\), \(\beta_2=-1\).
В этом случае выполняются равенства \(\zeta_1=-1\), \(\zeta_2=1\), \(\mathrm Z_+=\mathrm Z_-=1\).
Данные таблицы иллюстрируют утверждение \ref{par:3}.\ref{tm:3}.

\subsection
При получении вышеприведённого иллюстративного материала нами была использована
вычислительная методика, описанная в работе \cite{V}.


\begin{thebibliography}{9}
\bibitem[СФ]{VSh1} А.~А.~Владимиров, И.~А.~Шейпак. \emph{Самоподобные функции
в пространстве \(L_2[0,1]\) и задача Штурма--Лиувилля с сингулярным индефинитным
весом}// Матем.~сборник. "--- 2006. "--- Т.~197, \No~11. "--- С.~13--30.
\bibitem[ИШЛ]{VSh2} А.~А.~Владимиров, И.~А.~Шейпак. \emph{Индефинитная задача
Штурма--Лиувилля для некоторых классов самоподобных сингулярных весов}//
Труды~МИРАН им.~В.~А.~Стеклова. "--- 2006. "--- Т.~255. "--- С.~88--98.
\bibitem[SSM]{SV} M.~Solomyak, E.~Verbitsky. \emph{On a spectral problem related
to self-similar measures}// Bull.~London Math.~Soc. "--- 1995. "--- V.~27,
\No~3. "--- P.~242--248.
\bibitem[ВСЗ]{V} А.~А.~Владимиров. \emph{О вычислении собственных значений задачи
Штурма--Лиувилля с фрактальным индефинитным весом}// Журнал выч.~матем.
и матем.~физ. "--- 2007. "--- Т.~47, \No~8. "--- С.~1350--1355.
\end{thebibliography}
\end{document}